\documentclass[AMA,Times1COL]{WileyNJDv5} 

\articletype{Article Type}%

\received{Date Month Year}
\revised{Date Month Year}
\accepted{Date Month Year}
\journal{Journal}
\volume{00}
\copyyear{2023}
\startpage{1}

\begin{document}

\title{Unraveling Forward and Backward Source Problems for a Nonlocal Integrodifferential Equation: A Journey through Operational Calculus for Dzherbashian-Nersesian Operator}

\author[1]{Anwar Ahmad}

\author[2]{Muhammad Ali}

\author[1]{Salman A. Malik}

\authormark{AHMAD \textsc{et al.}}
\titlemark{OPERATIONAL CALCULUS FOR DZHERBASHIAN-NERSESIAN OPERATOR}

\address[1]{\orgdiv{Department of Mathematics}, \orgname{COMSATS University Islamabad}, \orgaddress{\street{Park Road, Tarlai Kalan}, \city{Islamabad 45550}, \country{Pakistan}}}

\address[2]{\orgdiv{Department of Sciences and Humanities}, \orgname{University of Computer and Emerging Sciences}, \orgaddress{\city{Islamabad}, \country{Pakistan}}}


\corres{Corresponding author Salman A. Malik, {\orgdiv{Department of Mathematics}, \orgname{COMSATS University Islamabad}, \orgaddress{\street{Park Road, Tarlai Kalan}, \city{Islamabad 45550}, \country{Pakistan}}} \email{salman\_amin@comsats.edu.pk}}



\abstract[Abstract]{This article primarily aims at introducing a novel operational calculus of Mikusi\'nski's type for the Dzherbashian-Nersesian operator.
	Using this calculus, we are able to derive exact solutions for the forward and backward source problems (BSPs) of a differential equation that features Dzherbashian-Nersesian operator in time and intertwined with nonlocal boundary conditions. The initial condition is expressed in terms of Riemann-Liouville integral (RLI). Solution is presented using Mittag-Leffler type functions (MLTFs). The outcomes related to the existence and uniqueness subject to certain conditions of regularity on the input data are established.}


\jnlcitation{\cname{%
\author{Taylor M.},
\author{Lauritzen P},
\author{Erath C}, and
\author{Mittal R}}.
\ctitle{On simplifying ‘incremental remap’-based transport schemes.} \cjournal{\it J Comput Phys.} \cvol{2021;00(00):1--18}.}

\maketitle

%

\section{Introduction and Problems Statement}\label{intro.mokin}
\numberwithin{equation}{section}
Merit for developing operational calculus is credited to the famous Polish mathematician Jan Mikusi\'nki. He presented the new approach for constructing an operational calculus for the differential operator during the 1950s \cite{Mikusinski}. This algebraic scheme revolves around the idea that the Laplace convolution can be interpreted as a multiplication within the continuous functions ring in the upper half-plane. This approach, known as the Mikusi\'nski operational calculus, has been effectively applied to various mathematical fields, such as the study of ordinary and partial differential equations. Many mathematicians (see, for example \cite{Ditkin-1957, Ditkin-Prudnikov, Dimovski,Meller}) extended the Mikusi\'nki scheme to create operational calculi for differential equations with variable coefficients, including the Bessel differential operator defined as follows
\begin{equation*}
	(\mathcal{B}g)(x):=x^{-\alpha}\prod^{n}_{i=1}\Big(\beta_{i}+\frac{1}{\alpha}x\frac{d}{dx}\Big)g(x).
\end{equation*}
In year 1994, Luchko and Srivastva developed an operational calculus for the Riemann-Liouville derivative (RLD), which allowed them to solve a Cauchy boundary value problem associated with a linear equation that involved RLDs \cite{Luchko-Srivastava-CAMWA}. The following year, Al-Bassam and Luchko constructed an operational calculus pertaining to the multiple Erdelyi-Kober operators \cite{Al-Bassam}. This operational calculus enabled them to solve the nonlocal integrodifferential equations with multiple Erdelyi-Kober operators. In 1999, Luchko and Gorenflo published an article in Acta Mathematica in which they introduced Mikusi\'nki's operational calculus for Caputo derivative \cite{Luchko-Gorenflo-1999}. 
In 2009, Hilfer and coauthors developed an operational calculus for the generalized RLD \cite{Hilfer-Luchko-Tomovski}. This scheme of operational calculus entails solving linear differential equations featuring constant coefficients that incorporate generalized RLD.

Despite being introduced in 1968 \cite{Dzhrbashyan-Nersesyan-1968}, the Dzherbashian-Nersesian operator has not received significant attention from the mathematical community.  However, the appearance of an English translation of Dzhrbashyan et al. \cite{Dzhrbashyan-Nersesyan-1968} in an article published in FCAA \cite{Dzhrbashyan-Nersesyan-2020} has sparked interest and discussion related to this operator. In addition to presenting several interesting results related to this operator, such as the construction of Laplace transform in Ahmad et al. \cite{Anwar-Ali-Malik-2021}, the authors have also explored the specific parameter choices that allow Dzherbashian-Nersesian operator to recover the renowned fractional operators like Riemann-Liouville, Caputo, and Hilfer derivatives. 

Forward and backward problems for FDEs arise in various areas of science. These kind of problems have also been studied for fractional diffusion equations in \cite{Kerbal-Kadirkulov-Kirane,Restrepo-Suragan,Ali-Malik-2017}. 
Recently, backward problems for FDEs involving Caputo and generalized RLDs have been extensively studied \cite{Hamidi-Kirane-Tfayli,Dib-Kirane,Al-Salti-Karimov}.
However, there is a limited body of literature pertaining to the forward and backward problems associated with FDEs that feature the Dzherbashian-Nersesian operator. To our knowledge, this marks the inaugural attempt to address a forward problem involving this operator, with only a handful of articles in the existing literature addressing backward problems for FDEs incorporating this operator. A few such instances include Berdyshev et al.\cite{Berdyshev-Kadirkulov-2016} and Karimov et al. \cite{Karimov-Bakhodirjon-2022}. Additionally, in \cite{Anwar-Ali-Malik-2021}, authors used the Fourier method to investigate backward problems involving source terms that are both space and time dependent.

Motivated by the persistent research activity in operational calculus, the renewed interest in Dzherbashian-Nersesian operator, and applications of backward problems, our study in this paper is of dual purpose. First of all, this chapter addresses the construction of operational calculus for Dzherbashian-Nersesian operator. Second, this calculus is used to address both forward and backward problems for FDE with Dzherbashian-Nersesian operator.

The subsequent sections of the article are structured in the following manner. Further in this section, we formulate the problems. Section \ref{Preliminaries.Operational} is dedicated to some basics definitions and results related to nonlocal integrodifferential operators and Mittag-Leffler function. Section \ref{sec: Operational} focuses the development of operational calculus for Dzherbashian-Nersesian operator. Both forward and backward problems for a FDE are addressed in Section \ref{sec: FSP and BSP}. Finally, the concluding remarks are given in Section \ref{sec: Conclusion}.
\subsection*{Formulation of Forward and Backward Problems}
This paper focuses on investing the following FDE:
\begin{equation}\label{probeq.oper}
	\mathcal{D}^{\varrho_{m}}_{t,0+}u(t,x)=u_{xx}(t,x)+F(t,x),\quad (t,x)\in (0,T)\times (0,1),
\end{equation}
contingent upon the subsequent boundary conditions
\begin{equation}\label{probbcs.oper}
	u(t,1)=0,\quad u_{x}(t,0)=u_{x}(t,1),\quad t\in (0,T),
\end{equation}
and initial condition given by
\begin{equation}\label{probic.oper}
	J^{1-\alpha_{1}}_{t,0+}u(t,x)\Big|_{t=0}=\phi(x),\quad  \alpha_{0}, x\in (0,1).
\end{equation}
The expression $\mathcal{D}^{\varrho_{m}}_{t,0+}$, such that $\varrho_{m} \in (0,1)$, denotes Dzherbashian-Nersesian operator (See Definition \ref{DherNer-def}).

The pursuit of finding $u(t,x)$, which adheres to the initial boundary value problem (IBVP) described by (\ref{probeq.oper})-(\ref{probic.oper}) while satisfying the conditions $t^{\alpha_{1}}u(t,.)\in C^{2}(0,1)(.,x)$ and $t^{\alpha_{1}}\mathcal{D}^{\varrho_{m}}_{t,0+}u(.,x) \in C(0,T)$, given the source term $F(t,x)$ and the initial data $\phi(x)$ are provided, is recognized as the establishment of a strong or classical solution to (\ref{probeq.oper})-(\ref{probic.oper}). This is commonly referred to as the forward problem.

We also explore the BSP within the context of IBVP (\ref{probeq.oper})-(\ref{probic.oper}). It is important to note that the source term $F(t,x)$ ,is considered to be solely dependent on spatial coordinate, denoted as $F(t,x):=f(x)$. Let us formally introduce the concept of solving BSP: We aim to identify the set $\big\{u(t,x), f(x)\big\}$ that satisfy the conditions, specifically, $t^{\alpha_{1}}u(t,.)$ $t^{\alpha_{1}}\mathcal{D}^{\varrho_{m}}_{t,0+} u(.,x)$, and $f$ belonging to the spaces $C^{2}[0,1]$, $C(0,T)$, and $\in C[0,1]$ respectively, for the problem (\ref{probeq.oper})-(\ref{probic.oper}). This solution of BSP presents challenges as it falls under the category of ill-posed problems according to Hadamard's definition. To uniquely determine the source term, an additional condition is required, commonly referred to as the overdetermined condition. We propose the following overdetermination condition, thereby enhancing the uniqueness of the source term determination:
\begin{equation}\label{probover.oper}
	u(T,x)=\psi(x).	
\end{equation}
\section{Preliminaries and Notations}\label{Preliminaries.Operational}
\numberwithin{equation}{section}
This section is dedicated to presenting the fundamental definitions and properties of some well-known nonlocal integrodifferential operators and Mittag-Leffler function. 
\begin{definition}\cite{Kilbas-Srivastava}, \cite{Samko-book}
	The RLI denoted as $J^{\alpha}_{t,0+}$, where $\alpha \in \mathbb{C}$ and $\mathcal{R}(\alpha)>0$, is given by
	\begin{align}\label{defRLI.oper}
		J^{\alpha}_{t,0+}g(t):=\frac{1}{\Gamma(\alpha)}\int_{0}^{t}\frac{g(\tau)}{(t-\tau)^{1-\alpha}}d\tau, \quad t>0.
	\end{align}
\end{definition}
The function denoted by $\Gamma(\cdot)$ refers to the second kind of Euler integral.
\begin{definition}\cite{Kilbas-Srivastava}, \cite{Samko-book} \label{RL-def}
	The RLD denoted as $D^{\alpha}_{t,0+}$, with order $\alpha \in (m-1,m)$ such that $\mathcal{R}(\alpha)\geqq 0$, is given by
	\begin{align*}
		D^{\alpha}_{t,0+}g(t):= \frac{d^{m}}{dt^{m}}J^{m-\alpha}_{t,0+}g(t), \quad m=\mathcal{R}(\alpha)+1,
	\end{align*}
	where $[\mathcal{R}(\alpha)]$ represents the integer part of $\mathcal{R}(\alpha)$.
\end{definition}
\begin{definition}\cite{Dzhrbashyan-Nersesyan-1968}\label{DherNer-def}
	The Dzherbashian-Nersesian operator $\mathcal{D}^{\varrho_{m}}_{t,0+}$ in relation to time with order $\varrho_{m} \in (0,m)$ defined as follows
	\begin{align}
		\mathcal{D}^{\varrho_{m}}_{t,0+} g(t):=J^{1-\alpha_{m}}_{t,0+}D^{\alpha_{m-1}}_{t,0+}D^{\alpha_{m-2}}_{t,0+}...D^{\alpha_{1}}_{t,0+}D^{\alpha_{0}}_{t,0+}g(t),\quad m\in \mathbb{Z^{+}},\quad t>0,\label{DzhrNer-mathematical}
	\end{align}
	where
	\begin{align*}
		\varrho_{m}=\sum_{j=0}^{m}\alpha_{j}-1>0,\quad \alpha_{j}\in (0,1].
	\end{align*}
\end{definition}
\begin{definition}\label{mittag}\cite{Mittag-Leffler-Sur}
	Defined by Magnus G\"{o}sta Mittag-Leffler, the classical Mittag-Leffler function can be expressed in the following manner:
	\begin{align*}
		E_{\alpha}(z):=\sum_{k=0}^{\infty}\frac{z^{k}}{\Gamma(\alpha k+ 1)},\quad \mathcal{R}(\alpha)>0,\; z\in \mathbb{C}.
	\end{align*}	
	Wiman provided the generalization of (\ref{mittag}) in \cite{Wiman} as follows:
	\begin{align*}
		E_{\alpha, \beta}(z):=\sum_{k=0}^{\infty}\frac{z^{k}}{\Gamma(\alpha k+ \beta)},\quad \mathcal{R}(\alpha),\:\mathcal{R}(\beta)>0, \:z \in\mathbb{C}.	
	\end{align*}
\end{definition}
\noindent Moreover, the MLTF is defined in the following manner
\begin{align*}
	e_{\alpha,\beta}(t,\lambda):=t^{\beta-1}E_{\alpha,\beta}(-\lambda t^{\alpha}),\quad \mathcal{R}(\alpha),\:\mathcal{R}(\beta),\;t,\;\lambda>0.
\end{align*}
\begin{lemma}\cite{Podlubny}\label{podlem} The relation given below is valid subject to certain conditions that $\alpha<2$, $\beta$ is any real number, $\mu$ is a value satisfying $\pi\alpha/2<\mu<min\{\pi,\pi\alpha\}$, $z$ is a complex number in a manner that it satisfies $|z|\geq0$,
	$\mu\leq |arg(z)|\leq\pi$ and $C_{1}$ is a real constant.
	\begin{align*}
		\big|E_{\alpha,\beta}(z)\big|\leq \frac{C_{1}}{1+|z|}.
	\end{align*}	
\end{lemma}	
\begin{lemma}\cite{Ali-Aziz-Salman-FCAA1}\label{AliMittagtypeLemma1}
	The MLTF $e_{\alpha,\alpha+1}(t,\lambda)$ exhibits the subsequent property:
	\begin{align*}
		e_{\alpha,\alpha+1}(t,\lambda)=\frac{1}{\lambda}\big(1-e_{\alpha,1}(t,\lambda)\big),\quad t, \lambda>0.
	\end{align*}
\end{lemma}
\begin{lemma}\cite{Ali-Aziz-Salman-FCAA1}\label{AliMittagtypeLemma1'}
	The subsequent property holds for the MLTF $e_{\alpha,1}(T,\lambda)$, for any $\alpha$ within the interval $(0,1)$:
	\begin{align*}
		\frac{1}{1-e_{\alpha,1}(T,\lambda)}\leq C_{2},\quad T,\: \lambda, \: C_{2}>0.
	\end{align*}
\end{lemma}	
\section{Operational calculus for Dzherbashian-Nersesian operator}\label{sec: Operational}
\numberwithin{equation}{section}
The section showcases the construction of operational calculus for the Dzherbashian-Nersesian operator. This calculus is afterwards in the paper plays a pivotal role in investigating (\ref{probeq.oper})-(\ref{probic.oper}) and (\ref{probeq.oper})-(\ref{probover.oper}). Furthermore, under certain fixation of parameters, the formula of operational calculus for Dzherbashian-Nersesian operator $\mathcal{D}^{\varrho_{m}}_{t,0+}$ interpolate those of fractional derivatives namely Riemann-Liouville, Caputo and Hilfer. 

However, we start with providing some important operational relations from the existing literature.
\begin{definition}\cite{Luchko-Srivastava-CAMWA}
	A function $g$ with real or complex values is considered to be a member of the space $C_{\gamma}$, where $\gamma$ is a real number, only if a real number $p$ greater than $\gamma$ exists in a manner that
	\begin{align*}
		g(z)=z^{p}g_{1}(z)	,\quad t>0,
	\end{align*}
	where $g_{1}\in C[0,\infty)$.
\end{definition}
The fact that $C_{\gamma}$ forms a vector space is evident. Moreover, the collection of $C_{\gamma}$ is arranged in an order based on inclusion given by the subsequent relation
\begin{align*}
	C_{\gamma}\subset C_{\delta},\quad \iff \gamma \ge \delta.
\end{align*}
\begin{theorem}\cite{Luchko-Srivastava-CAMWA}
	The RLI $J^{\alpha}_{t,0+}$, such that $\alpha>0$, is linear mapping of the space $C_{\gamma}$, $\gamma \geq -1$, into itself, i.e.,  
	\begin{align*}
		J^{\alpha}_{t,0+}:C_{\gamma}\rightarrow C_{\gamma+\alpha}\subset C_{\gamma}.
	\end{align*}
\end{theorem} 
It is a well established fact that $J^{\alpha}_{t,0+}$, $\alpha>0$, possesses the subsequent convolution representation $C_{\gamma}$, $\gamma \geq -1$:
\begin{align*}
	(J^{\alpha}_{t,0+}y)(z)=(h_{\alpha}* y)(z), \quad h_{\alpha}(z):=\frac{z^{\alpha}-1}{\Gamma(\alpha)}, \quad y\in C_{\gamma}.
\end{align*}
where $*$ denotes the Laplace convolution defined as under:
\begin{align*}
	(g* f)(z)=\int_{0}^{z}g(z-t)f(t)dt,\quad z>0,	
\end{align*}
\begin{definition}\cite{Hilfer-Luchko-Tomovski}
	If $D^{\alpha}_{t,0+}y\in C_{-1}$ for all $\alpha$ between $n-1$ and $n$, then $y$ belongs to space $\Omega^{\alpha}_{-1}$ and has $\alpha\geq 0$.
\end{definition}
\begin{theorem}\cite{Luchko-Srivastava-CAMWA}
	The commutative ring $(C_{-1}, *, +)$ formed by the Laplace convolution $*$ and ordinary addition in the space $C_{-1}$ is free from zero divisors. 
\end{theorem}
By applying traditional Mikusi\'nski operational calculus, one can extend this ring to the field of convolution quotients denoted as  $M_{-1}$:
\begin{align*}
	M_{-1}:=C_{-1}\times (C_{-1}\cup \{0\})/\sim,
\end{align*}
with the relation of equivalence $(\sim)$ described as
\begin{align*}
	(f,g)\sim (f_{1},g_{1})	\equiv (f* g_{1})(t)=(g* f_{1})(t).
\end{align*}
The elements within the field $M_{1}$ can be viewed as convolution quotients $f/g$, for the sake of simplicity. However, the standard operations of addition and multiplication continue to hold their definitions in $M_{1}$ in the following manner:
\begin{align}
	\frac{f}{g}+\frac{f_{1}}{g_{1}}:=\frac{f* g_{1}+g* f_{1}}{g* g_{1}},\label{add}\\
	\frac{f}{g}\circ\frac{f_{1}}{g_{1}}:=\frac{f* f_{1}}{g* g_{1}}.\label{multi}
\end{align}
\begin{theorem}\cite{Luchko-Srivastava-CAMWA}
	The commutative field $(M_{-1},.,+)$ can be constructed by defining the addition (\ref{add}) and multiplication (\ref{multi}) operations within $M_{1}$. Using a mapping that assigns $(\alpha > 0)$, we can embed the ring $C_{-1}$ in $M_{-1}$ as follows:
	\begin{align*}
		f \longmapsto \frac{h_{\alpha}* f}{h_{\alpha}}.
	\end{align*}
\end{theorem}
It is possible to define the scalar multiplication operation in $M_{-1}$ using the relation given below where the scalar $\lambda$ is taken from $\mathbb{R}$ or $\mathbb{C}$
\begin{align*}
	\lambda\frac{f}{g}:=\frac{\lambda f}{g},\quad \frac{f}{g}\in M_{-1}.
\end{align*}
In $M_{-1}$, the definition of the operation involving a constant function and multiplication is as follows:
\begin{align*}
	\{\lambda\}	\circ \frac{f}{g}=\frac{\lambda h_{\alpha+1}}{h_{\alpha}}* \frac{f}{g}=\{1\}* \frac{f}{g}.
\end{align*}
Specifically, we assume an element in $M_{-1}$ denoted by $I$. This element can be defined as the ratio of $h_{\alpha}$ to itself and serves as the element of identity in $M_{-1}$ when using the multiplication operation given as follows:
\begin{align*}
	I \circ \frac{f}{g}=\frac{h_{\alpha}* f}{h_{\alpha}* g}=\frac{f}{g}.
\end{align*}
In the conventional theory of generalized functions, the Dirac $\delta$-function is analogous to the element of identity $I$ in $M_{-1}$, as illustrated by the last formula.
\begin{definition}\cite{Luchko-Srivastava-CAMWA}
	In $M_{-1}$, the reciprocal to the element $h_{\alpha}$ is denoted as $S_{\alpha}$. It is recognized as the algebraic inverse of the RLI operator $J^{\alpha}_{t,0+}$, i.e.,
	\begin{align*}
		S_{\alpha}=\frac{I}{h_{\alpha}}\equiv \frac{h_{\alpha}}{h_{\alpha}* h_{\alpha}}\equiv \frac{h_{\alpha}}{2h_{\alpha}},
	\end{align*}
	where $I$ represents the element of identity in $M_{-1}$ under the operation of multiplication, given by the ratio of $h_{\alpha}$ to itself.
\end{definition}
The representation of RLI $J^{\alpha}_{t,0+}$ involves multiplication (convolution) with the function $h_{\alpha}$, as mentioned in \cite{Dimovski}. This multiplication (convolution) takes place in the ring $C_{-1}$. By the virtue of embedding of $C_{-1}$ in $M_{-1}$ of convolution quotients, it can be expressed in alternative form as:
\begin{align*}
	(J^{\alpha}_{t,0+})y(x)=\frac{1}{S_{\alpha}}\circ y.
\end{align*}
The operational calculus for the Dzherbashian-Nersesian operator $\mathcal{D}^{\varrho_{m}}_{t,0+}$ is significantly aided by the following Theorem in solving the nonlocal integrodifferential equations:
\begin{theorem}\cite{Luchko-2020-FCAA, Luchko-2020-Mathematics} 
Let $y\in \Omega^{\varrho_{m}}_{-1}$, $0<\varrho_{m}\leq m$, $m \in \mathbb{Z^{+}}$. Then the RLI (\ref{defRLI.oper}) and the Dzherbashian-Nersesian operator (\ref{DzhrNer-mathematical}) are connected with each other by the following relation
\begin{align}\label{RLIDzh.oper}
	(J^{\varrho_{m}}_{t,0+}\mathcal{D}^{\varrho_{m}}_{t,0+}y)(x)=y(x)-\sum_{k=0}^{m-1}\frac{x^{\varrho_{k}}}{\Gamma(\varrho_{k}+1)}(\mathcal{D}^{\varrho_{k}}_{t,0+}y)(x)\Big|_{x=0}, \quad x>0.
\end{align}
\end{theorem}
\begin{theorem}\label{inverse.oper}
	Let function $y \in \Omega^{\varrho_{m}}_{-1}$, where $0<\varrho_{m}\leq m$, $m\in \mathbb{Z^{+}}$. In this case, the Dzherbashian-Nersesian operator $\mathcal{D}^{\varrho_{m}}_{t,0+}$ can be expressed as multiplication in $M_{-1}$ of convolution quotients as follows:
	\begin{align}\label{dzhrmain.oper}
		(\mathcal{D}^{\varrho_{m}}_{t,0+}y)(x)=S_{\varrho_{m}}\circ y(x)-\sum_{k=0}^{m-1}S_{\varrho_{m}-\varrho_{m-k}-1}\circ (\mathcal{D}^{\varrho_{k}}_{t,0+}y)(x)\Big|_{x=0}, \quad x>0.
	\end{align}
\end{theorem}  
\begin{proof}
	To obtain the formula (\ref{dzhrmain.oper}), initially, the embedding of the ring $C_{-1}$ in $M_{-1}$ is carried out. Subsequently, the relationship (\ref{RLIDzh.oper}) undergoes multiplication with the algebraic inverse of the RLI, denoted by $S_{\alpha}$. This gives the desired relation without any differences from formula (\ref{dzhrmain.oper}).
\end{proof}
\begin{remark} 
	On setting the parameters $\alpha_{0}=1+\alpha-m$ and $\alpha_{1}=\alpha_{2}=...=\alpha_{k}=1$ in Equation (\ref{dzhrmain.oper}), we obtain the relation of operational calculus for Riemann-Liouville derivative \cite{Luchko-Srivastava-CAMWA}, i.e.,
	\begin{align*}
		\mathcal{D}^{\varrho_{m}}_{t,0+} y=S_{\alpha}\circ y-S_{\alpha}\circ y_{\alpha},
	\end{align*}
	where $$y_{\alpha}(x):=\sum_{k=0}^{m-1}\frac{x^{k+\alpha-m}}{\Gamma(k+\alpha-m+1)}\lim_{x\rightarrow 0+}D^{k+\alpha-m}_{t,0+}y(x).$$
\end{remark}
\begin{remark}
	The Equation (\ref{dzhrmain.oper}) admits the relation of operational calculus for Caputo derivative \cite{Luchko-Gorenflo-1999}if we take parameters $\alpha_{1}=\alpha_{2}=...=\alpha_{k-1}=1$ and $\alpha_{k}=1+\alpha-k$, i.e.,
	\begin{align*}
		\mathcal{D}^{\varrho_{m}}_{t,0+}y=S_{\alpha}\circ y-S_{\alpha}\circ y_{\alpha},
	\end{align*}
	where $$y_{\alpha}(x):=\sum_{k=0}^{m-1}\frac{x^{k}}{\Gamma(k+1)}\lim_{x\rightarrow 0+}y^{(k)}(x).$$
\end{remark}
\begin{remark}
	The relation for operational calculus for Hilfer derivative \cite{Hilfer-Luchko-Tomovski} is interpolated if we fix the parameters $\alpha_{0}=1-(m-\alpha)(1-\beta)$ and $\alpha_{1}=\alpha_{2}=...=\alpha_{k}=1$ in Equation (\ref{dzhrmain.oper}), i.e.,
	\begin{align*}
		\mathcal{D}^{\varrho_{m}}_{t,0+}y=S_{\alpha}\circ y-S_{\alpha}\circ y_{\alpha,\beta},
	\end{align*}
	where $$y_{\alpha,\beta}(x):=\sum_{k=0}^{m-1}\frac{x^{k-m+\alpha-\beta\alpha+\beta m}}{\Gamma(k-m+\alpha-\beta\alpha+\beta m+1)}\lim_{x\rightarrow 0+}\frac{d^{k}}{dx^{k}}(J^{(1-\beta)(m-\alpha)}_{t,0+}y)(x).$$
\end{remark}
\begin{lemma}\cite{Luchko-Gorenflo-1999}\label{Luchko-Caputo-Lemma}
For $\alpha_{i}$, $m_{i}$, $\sigma$, $\eta$ $>0$, we have
\begin{align*}
\frac{S_{-\eta\sigma}}{1-\sum_{k=1}^{\infty}}=z^{\eta\sigma-1}E_{\alpha_{1}\sigma,...,\alpha_{n}\sigma},\eta\sigma(m_{1}z^{\alpha_{1}\eta},...,m_{n}z^{\alpha_{n}\eta}).
\end{align*}
\end{lemma}
\section{Forward and Backward Source Problems}\label{sec: FSP and BSP}
In this section, we will utilize the operational calculus derived in the previous section to examine both the forward and backward problems (\ref{probeq.oper})-(\ref{probic.oper}) and (\ref{probeq.oper})-(\ref{probover.oper}), respectively. However, we begin our discussion by exploring the spectral problem associated with (\ref{probeq.oper}) and (\ref{probbcs.oper}), as it plays a crucial role in addressing the problems delineated by (\ref{probeq.oper})-(\ref{probic.oper}) and (\ref{probeq.oper})-(\ref{probover.oper}).
\subsection{Spectral Problem}\label{spectral.oper}
\numberwithin{equation}{section}
Our approach to solving the BSP (\ref{probeq.oper})-(\ref{probover.oper}) involves the Fourier method, also referred to as separation of variables. The nonlocal boundary conditions in this problem render the spectral non-self-adjoint.

Below is the spectral problem for (\ref{probeq.oper}) contingent upon boundary conditions (\ref{probbcs.oper})
\begin{equation}\label{probspec.oper}
	\begin{cases}
		X^{''}(x)=-\lambda X(x),	\quad x\in (0,1),\\	
		X(1)=0,\quad X'(0)=X'(1).
	\end{cases}	
\end{equation}	
The set of eigenfunctions denoted as $\{X_{0}, X_{1k}\}$ for (\ref{probspec.oper}) with respective eigenvalues $\lambda_{0}=0$, $\lambda_{k}=(2\pi k)^{2}$, along with the associated eigenfunction $ X_{2k}$ (see \cite{Moiseev-1999}), can be expressed as follows
\begin{align}\label{eigenfpec.oper}
	X_{0}(x)=2(1-x),\quad X_{1k}(x)	=4(1-x)\cos (2\pi kx),\quad X_{2k}(x)=4\sin (2\pi kx).
\end{align}
Under the spectral problem (\ref{probspec.oper}), the adjoint problem is presented as follows
\begin{equation}\label{probadj.oper}
	\begin{cases}
		Y^{''}(x)=-\lambda Y(x),	\quad x\in (0,1),\\	
		Y(0)=0,\quad Y'(0)=Y'(1).
	\end{cases}	
\end{equation}	
The set of eigenfunctions denoted as $\{Y_{0}, Y_{1k}\}$ for (\ref{probadj.oper}) with respective eigenvalues $\lambda_{0}=0$, $\lambda_{k}=(2\pi k)^{2}$, together with the associated eigenfunction $ Y_{2k}$ (see \cite{Moiseev-1999}), can be expressed as follows
\begin{align}\label{eigenfadj.oper}
	Y_{0}(x)=1,\quad 	Y_{1k}(x)=\cos(2\pi kx), \quad Y_{2k}(x)=x\sin(2\pi kx).
\end{align}
\begin{lemma}\label{Riesz.oper}\cite{Ionkin-Moiseev}
	The systems of functions illustrated in (\ref{eigenfpec.oper}) and (\ref{eigenfadj.oper}) establish Riesz basis in $L^{2}(0,1)$.
\end{lemma}
\subsection{Forward Problem}
This subsection is dedicated to the study of forward problem described by (\ref{probeq.oper})-(\ref{probic.oper}). We demonstrate, subject to specific assumptions given in Theorem \ref{fp.oper}, the existence of a unique classical solution to the forward problem. We will outline the outcomes related to the existence and uniqueness of solutions for the forward problem.

The forward problem (\ref{probeq.oper})-(\ref{probic.oper}) is linear. This allows us to express the solution using the following formulation:
\begin{align*}
	u(t,x)=v(t,x)+w(t,x),
\end{align*}
where $v(t,x)$ is the solution of following associated homogeneous problem, i.e.,
\begin{equation}\label{probeq.homo.oper}
	\mathcal{D}^{\varrho_{m}}_{t,0+}v(t,x)=v_{xx}(t,x),\quad (t,x)\in (0,T)\times (0,1),
\end{equation}
\begin{equation}\label{probbcs.homo.oper}
	v(t,1)=0,\quad v_{x}(t,0)=v_{x}(t,1),\quad t\in (0,T),
\end{equation}
\begin{equation}\label{probic.homo.oper}
	J^{1-\alpha_{1}}_{t,0+}v(t,x)\Big|_{t=0}=\phi(x),\quad  \alpha_{1}, x\in (0,1),
\end{equation}
and $w(t,x)$ is the solution of following associated non-homogeneous problem, i.e.,
\begin{equation}\label{probeq.nonhomo.oper}
	\mathcal{D}^{\varrho_{m}}_{t,0+}w(t,x)=w_{xx}(t,x)+f(x),\quad (t,x)\in (0,T)\times (0,1),
\end{equation}
\begin{equation}\label{probbcs.nonhomo.oper}
	w(t,1)=0,\quad w_{x}(t,0)=u_{x}(t,1),\quad t\in (0,T),
\end{equation}
\begin{equation}\label{probic.nonhomo.oper}
	J^{1-\alpha_{1}}_{t,0+}w(t,x)\Big|_{t=0}=0,\quad  \alpha_{1}, x\in (0,1).
\end{equation}
\begin{theorem}\label{fp.oper}
	For $\varrho_{m}\in (0,1)$ and $\phi(x) \in C^{2}(0,1)$  in a manner that it satisfies $\phi(1)=0$, $\phi'(0)=\phi'(1)$, the forward problem (\ref{probeq.homo.oper})-(\ref{probic.homo.oper}) possesses a unique classical solution given by:
	\begin{align}
		v(t,x)=&\varphi_{0}\frac{t^{\alpha_{0}-1}}{\Gamma(\alpha_{0})}	2(1-x)+\sum_{k=1}^{\infty} \Big[\Big\{\varphi_{1k} e_{\alpha_{0}+\alpha_{1}-1,\alpha_{0}}(t,\lambda_{k})\Big\}4(1-x)\cos (2\pi kx)\notag\\&+\Big\{e_{\alpha_{0}+\alpha_{1}-1,\alpha_{0}}(t,\lambda_{k})\varphi_{2k}+2\sqrt{\lambda_{k}}e_{\alpha_{0}+\alpha_{1}-1,\alpha_{0}+\alpha_{1}-1}(t,\lambda_{k})*e_{\alpha_{0}+\alpha_{1}-1,\alpha_{0}}(t,\lambda_{k})\varphi_{1k}\Big\}4\sin (2\pi kx)\Big],
	\end{align}
	where $\varphi_{i}:=\langle \varphi(x),Y_{i}(x)\rangle, \quad i \in \mathbb{Z}^{+}\cup \{0\}.$
\end{theorem} 
\begin{proof}
	The solution to the (\ref{probeq.homo.oper})-(\ref{probic.homo.oper}) is succeeded by the establishment of the existence and uniqueness of the obtained solution.
	\subsubsection*{Construction of the Solution:}
	Due to the Lemma \ref{Riesz.oper}, the set $\big\{X_{0}(x), X_{1k}(x), X_{2k}(x)\big\}$ forms a Riesz basis within the space $L^{2}(0,1)$, we are enabled to present the expression for $v(t,x)$ as follows:
	\begin{align}
		v(t,x)=v_{0}(t)X_{0}(x)+\sum_{k=1}^{\infty} \big(v_{1k}(t)X_{1k}(x)+v_{2k}(t)X_{2k}(x)\big),\label{u.homo.oper}
	\end{align}
	Utilizing the Equation (\ref{u.homo.oper}) within the framework of the Equation (\ref{probeq.homo.oper}), and considering the fact the bi-orthogonality of the function of sets $\big\{X_{0}(x), X_{1k}(x), X_{2k}(x)\big\}$ and $\big\{Y_{0}(x), Y_{1k}(x), Y_{2k}(x)\big\}$ in $L^{2}(0,1)$, we proceed to derive the subsequent system of FDEs:
	\begin{align}
		\mathcal{D}^{\varrho_{m}}_{t,0+}v_{0}(t)&=0,\label{fdeu0.homo.oper}\\
		\mathcal{D}^{\varrho_{m}}_{t,0+}v_{1k}(t)+\lambda_{k}v_{1k}(t)&=0,\label{fdeu1k.homo.oper}\\
		\mathcal{D}^{\varrho_{m}}_{t,0+}v_{2k}(t)-2\sqrt{\lambda_{k}}v_{1k}(t)+\lambda_{k}v_{2k}(t)&=0.\label{fdeu2k.homo.oper}
	\end{align}
	Using Theorem \ref{inverse.oper} in (\ref{fdeu0.homo.oper})-(\ref{fdeu2k.homo.oper}), we have
	\begin{align}
		S_{\alpha_{0}+\alpha_{1}-1}\circ v_{0}(t)&=S_{\alpha_{1}}\circ \phi_{0},\label{invuo.homo.oper}\\
		(S_{\alpha_{0}+\alpha_{1}-1}+\lambda_{k})\circ v_{1k}(t)&=S_{\alpha_{1}}\circ \phi_{1k},\label{invu1k.homo.oper}\\
		(S_{\alpha_{0}+\alpha_{1}-1}+\lambda_{k})\circ v_{2k}(t)-2\sqrt{\lambda_{k}}\circ v_{1k}(t)&=S_{\alpha_{1}}\circ \phi_{2k}.\label{invu2k.homo.oper}
	\end{align}
	Using Lemma \ref{Luchko-Caputo-Lemma} in (\ref{invuo.homo.oper})-(\ref{invu2k.homo.oper}), we have
	\begin{align}
		v_{0}(t)=&\frac{t^{\alpha_{0}-1}}{\Gamma(\alpha_{0})}\varphi_{0},\label{spaceu0.homo.oper}\\
		v_{1k}(t)=&e_{\alpha_{0}+\alpha_{1}-1,\alpha_{0}}(t,\lambda_{k})\varphi_{1k},\label{spaceu1k.homo.oper}\\
		v_{2k}(t)=&e_{\alpha_{0}+\alpha_{1}-1,\alpha_{0}}(t,\lambda_{k})\varphi_{2k}+2\sqrt{\lambda_{k}}e_{\alpha_{0}+\alpha_{1}-1,\alpha_{0}+\alpha_{1}-1}(t,\lambda_{k})*v_{1k}(t).\label{spaceu2k.homo.oper}
	\end{align}
	Prior to discussing the existence of the solution for the IBVP described by Equations (\ref{probeq.homo.oper})-(\ref{probic.homo.oper}), it is pertinent to state the following Lemma:
	\begin{lemma}\label{estlemma.oper}
		For $g\in C^{n}([0,1])$ in a manner that $g(1)=0$, $g^{(1)}(0)=g^{(1)}(1)$, $g(1)=0$ $g^{(3)}(0)=g^{(3)}(1)$ and we have the following relations:
		\begin{align*}
			\big|g_{1k}\big|\leq \frac{1}{k^{n}}\|g^{(n)}\|, \quad \big|g_{2k}\big|\leq \frac{n+1}{k^{n}}\|g^{(n)}\|,\quad n=1,2,3,4.
		\end{align*}
	\end{lemma} 
	Representation $\|.\|$ denotes the norm in $L^{2}(0,1)$ which is defined by
	\begin{align*}
		\|.\|:=\sqrt{\langle .,. \rangle}
	\end{align*}
	where $\langle .,. \rangle$ denotes the inner product which is given by $\langle f,g \rangle:=\int_{0}^{1}f(x)g(x)dx$.
	
	\subsubsection*{Existence of the Solution:}
	To prove the existence of the classical solution, we will prove that the series corresponding to $t^{\alpha_{1}}u(t,x)$, $t^{\alpha_{1}}\mathcal{D}^{\varrho_{m}}u(t,x)$, and $t^{\alpha_{1}}u_{xx}(t,x)$ represent continuous functions.
	
	\noindent Taking into account the CBSI, (\ref{spaceu0.homo.oper}) and the fact that $Y_{0}(x) = 1$, we obtain
	\begin{align}
		t^{\alpha_{1}}\big|v_{0}(t)\big|\leq \frac{t^{\alpha_{0}+\alpha_{1}-1}\|\phi\|}{\Gamma(\alpha_{0})}.
	\end{align}
	On using Lemma \ref{podlem}, Equation (\ref{spaceu1k.homo.oper}), we have
	\begin{align}
		\big|v_{1k}(t)\big|\leq 	\frac{C_{1}\big|\phi_{1k}\big|}{\lambda_{k}t^{\alpha_{1}}}.
	\end{align}
	Using the Lemma \ref{estlemma.oper}, CBSI and the property $\big|\frac{t}{T}\big|\leq 1$, we have
	\begin{align}\label{spaceu1k'.homo.oper}
		t^{\alpha_{1}}\big|v_{1k}(t)\big|\leq 	\frac{C_{1}\|\phi\|}{k^{2}}.
	\end{align}
	Likewise due to the Lemma \ref{podlem} along with Equations (\ref{spaceu2k.homo.oper}) and (\ref{spaceu1k'.homo.oper}), we have
	\begin{align}
		t^{\alpha_{1}}\big|v_{2k}(t)\big|\leq	\frac{C_{1}\|\phi\|}{k^{2}}+\frac{2C^{2}_{1}\|\phi\|}{k^{3}}.\label{spaceu2k'.homo.oper}
	\end{align}
	We still need to demonstrate the uniform convergence of the series associated with $\mathcal{D}^{\varrho_{m}}_{t,0+}v(t,x)$ within the interval $[\varepsilon, T]$. To establish this we employ the following lemma from \cite{Ali-Aziz-Salman-FCAA1}:
	
	\begin{lemma}\label{DzhrSamkolemma}
		Suppose that for each $i\in\mathbb{Z^{+}}$, we have a set $g_{i}$ of functions defined on $(0,b]$, that adhere to the subsequent conditions: 
		\begin{enumerate}
			\item derivatives $D^{\alpha_{0}}_{t,0+}g_{i}(t)$,
			$D^{\alpha_{1}}_{t,0+}D^{\alpha_{0}}_{t,0+}g_{i}(t)$,...,
			$D^{\alpha_{m-1}}_{t,0+}...D^{\alpha_{0}}_{t,0+}g_{i}(t)$ for $i\in \mathbb{Z^{+}}, t\in (0,b]$ exist,
			\item the series $\sum_{i=1}^{\infty}g_{i}(t)$ and $\sum_{i=1}^{\infty}D^{\alpha_{0}}_{t,0+}g_{i}(t)$, $\sum_{i=1}^{\infty}D^{\alpha_{1}}_{t,0+}D^{\alpha_{0}}_{t,0+}g_{i}(t)$,...,
			$\sum_{i=1}^{\infty}D^{\alpha_{m-1}}_{t,0+}... D^{\alpha_{1}}_{t,0+}D^{\alpha_{0}}_{t,0+}g_{i}(t)$ are
			uniformly convergent on the interval $[a+\varepsilon,b]$ for any $\varepsilon>0$.
		\end{enumerate}
		Then
		\begin{align*}
			\mathcal{D}^{\varrho_{m}}_{t,0+}	 \displaystyle\sum_{i=1}^{\infty}g_{i}(t)=\displaystyle\sum_{i=1}^{\infty}\mathcal{D}^{\varrho_{m}}_{t,0+}g_{i}(t).
		\end{align*}
	\end{lemma}
	Therefore, to ensure the continuity of $\mathcal{D}^{\varrho_{m}}_{t,0+}v(t,x)$, it is necessary for the series representations of $v(t,x)$ and $D^{\alpha_{0}}v(t,x)$ to exhibit uniform convergence. In view of (\ref{spaceu1k'.homo.oper}) and (\ref{spaceu2k'.homo.oper}), $v(t,x)$ is already continuous. Therefore, it suffices to prove the continuity of $D^{\alpha_{0}}_{t,0+}v(t,x)$.
	
	One can observe that $D^{\alpha_{0}}_{t,0+}v(t,x)$ is bounded from above by a series that converges. Therefore, using the WMT, it can be concluded that $D^{\alpha}_{t,0+}v(t,x)$ denotes a function which is continuous.
	
	Moreover, from Equations (\ref{fdeu0.homo.oper})-(\ref{fdeu2k.homo.oper}), we have the following estimates
	\begin{align}
		t^{\alpha_{1}}\big|\mathcal{D}^{\varrho_{m}}_{t,0+}v_{0}(t)\big|\leq& \frac{\Gamma(\alpha_{0}+\alpha_{1})\|\phi\|}{\Gamma(\alpha_{0})},\label{estsumu0.homo.oper}\\
		t^{\alpha_{1}}\big|\mathcal{D}^{\varrho_{m}}_{t,0+}v_{1k}(t)\big|\leq& \sum_{k=1}^{\infty}\frac{C_{1}\|\phi^{(2)}\|}{k^{2}},\label{estsumu1k.homo.oper}\\
		t^{\alpha_{1}}\big|\mathcal{D}^{\varrho_{m}}_{t,0+}v_{2k}(t)\big|\leq& \sum_{k=1}^{\infty}\Big(\frac{C_{1}\|\phi^{(2)}\|}{k^{2}}+\frac{2C_{1}\|\phi^{(2)}\|}{k^{3}}\Big).\label{estsumu2k.homo.oper}
	\end{align}
	By virtue of the Lemma \ref{DzhrSamkolemma} and (\ref{estsumu0.homo.oper})-(\ref{estsumu2k.homo.oper}), we notice that $t^{\alpha_{1}}\mathcal{D}^{\varrho_{m}}_{t,0+}v(t,x)$ has a finite upper bound. Therefore, by WMT $t^{\alpha_{1}}\mathcal{D}^{\varrho_{m}}_{t,0+}v(t,x)$ converges.
	
	\noindent Similarly for $v_{xx}(t,x)$, we have
	\begin{align*}
		t^{\alpha_{1}}\big|v_{xx}(t,x)\big|\leq 48\sum_{k=1}^{\infty}\frac{C_{1}\|\phi\|}{k^{2}},
	\end{align*}
	which on employing WMT denotes a continuous function.
	\subsubsection*{Uniqueness of the Solution:}
	Suppose $\bar{v}(t,x)=v_{1}(t,x)-v_{2}(t,x)$, where $v_{1}(t,x)$ and $v_{2}(t,x)$ represent two solutions of IBVP defined by (\ref{probeq.homo.oper})-(\ref{probic.homo.oper}). The function $v(t,x)$ adheres to the following system:
	\begin{equation}\label{probeq.homo.oper.uniq}
		\mathcal{D}^{\varrho_{m}}_{t,0+}\bar{v}(t,x)=\bar{v}_{xx}(t,x),\quad (t,x)\in (0,T)\times (0,1),
	\end{equation}
	\begin{equation}\label{probbcs.homo.oper.uniq}
		\bar{v}(t,1)=0,\quad \bar{v}_{x}(t,0)=\bar{v}_{x}(t,1),\quad t\in (0,T),
	\end{equation}
	\begin{equation}\label{probic.homo.oper.uniq}
		J^{1-\alpha_{1}}_{t,0+}\bar{v}(t,x)\Big|_{t=0}=0,\quad  \alpha_{1}, x\in (0,1),
	\end{equation}
	Taking into account the bi-orthogonality property within sets of functions $\{X_{0}(x), X_{1k}(x), X_{2k}(x)\}$ and $\{Y_{0}(x), Y_{1k}(x), Y_{2k}(x)\}$ in $L^{2}(0,1)$, we proceed to formulate the following system of FDEs:
	\begin{align}
		\mathcal{D}^{\varrho_{m}}_{t,0+}\bar{v}_{0}(t)&=0,\label{fdeu0.homo.oper.uniq}\\
		\mathcal{D}^{\varrho_{m}}_{t,0+}\bar{v}_{1k}(t)+\lambda_{k}\bar{v}_{1k}(t)&=0,\label{fdeu1k.homo.oper.uniq}\\
		\mathcal{D}^{\varrho_{m}}_{t,0+}\bar{v}_{2k}(t)-2\sqrt{\lambda_{k}}\bar{v}_{1k}(t)+\lambda_{k}\bar{v}_{2k}(t)&=0.\label{fdeu2k.homo.oper.uniq}
	\end{align}
	In the view of (\ref{probeq.homo.oper.uniq}), we obtain
	$\bar{v}_{0}(t)=0,\quad \bar{v}_{1k}(t)=0, \quad \bar{v}_{2k}(t)=0$.
	
	As a result, we have
	
	$$\bar{v}_{0}(t)=0$$.
	
	The solution of the problem outlined in (\ref{probeq.nonhomo.oper})--(\ref{probic.nonhomo.oper}) can be obtained through the application of Duhamel's principle as referenced in \cite{Umarov-2012,Umarov-Saidamatov-2007} and is given by:
	\begin{align*}
		w(t,x)=\int_{0}^{t} u^{w}(t,\tau,x)d\tau,
	\end{align*}
	where $u^{w}(t,\tau,x)$ is the solution of the subsequent system
	\begin{equation}\label{probeq.nonhomo.duhamel.oper}
		\mathcal{D}^{\varrho_{m}}_{t,0+}u^{w}(t,\tau,x)=u^{w}_{xx}(t,\tau,x),\quad (t,x)\in (0,T)\times (0,1),
	\end{equation}
	\begin{equation}\label{probbcs.nonhomo.duhamel.oper}
		u^{w}(t,\tau,1)=0,\quad u^{w}_{x}(t,\tau,0)=u^{w}_{x}(t,\tau,1),\quad t\in (0,T),
	\end{equation}
	\begin{equation}\label{probic.nonhomo.duhamel.oper}
		J^{1-\alpha_{1}}_{t,0+}u^{w}(t,\tau,x)\Big|_{t=\tau}=D^{2-\alpha_{0}-\alpha_{1}}_{t,0+}f(\tau,x),\quad  \alpha_{0}, \alpha_{1}, x\in (0,1),
	\end{equation}
	Applying the similar approach, we derive the solution for initial boundary value problem described in (\ref{probeq.nonhomo.duhamel.oper})-(\ref{probic.nonhomo.duhamel.oper}), resulting in:
	\begin{align*}
		w(t,x)=&\frac{t^{\alpha_{0}-1}}{\Gamma(\alpha_{0})}\varphi_{0}X_{0}(x)D^{2-\alpha_{0}-\alpha_{1}}_{t,0+}f_{0}(\tau)+\sum_{k=1}^{\infty} \Big(e_{\alpha_{0}+\alpha_{1}-1,\alpha_{0}}(t,\lambda_{k})\varphi_{1k}X_{1k}(x)D^{2-\alpha_{0}-\alpha_{1}}_{t,0+}f_{1k}(\tau)\\
		&+\big(e_{\alpha_{0}+\alpha_{1}-1,\alpha_{0}}(t,\lambda_{k})\varphi_{2k}D^{2-\alpha_{0}-\alpha_{1}}_{t,0+}f_{2k}(\tau)+2\sqrt{\lambda_{k}}e_{\alpha_{0}+\alpha_{1}-1,\alpha_{0}+\alpha_{1}-1}(t,\lambda_{k})*v_{1k}(t)D^{2-\alpha_{0}-\alpha_{1}}_{t,0+}f_{2k}(\tau)\big)X_{2k}(x)\Big).
	\end{align*}
\end{proof}
\subsection{Backward Problem}
Operational calculus constructed in the previous subsection will be used to solve the space dependent BSP (\ref{probeq.oper})-(\ref{probover.oper}). We begin by establishing the theorem addressing the BSP represented by (\ref{probeq.nonhomo.oper})-(\ref{probover.oper}).
\smallskip

\begin{theorem}
	For $\varrho_{m}\in (0,1)$ and \\
	$1$. $\phi(x) \in C^{2}(0,1)$  in a manner that it satisfies $\phi(1)=0$, $\phi'(0)=\phi'(1)$,\\
	$2$. $\psi(x) \in C^{4}(0,1)$  in a manner that it satisfies $\psi(1)=0$, $\psi'(0)=\psi'(1)$, $\psi''(1)=0$, $\psi'''(0)=0=\psi'''(1)$,\\
	the backward problem (\ref{probeq.oper})-(\ref{probover.oper}) possesses a unique regular solution given by:
\end{theorem} 
\begin{proof}
	The solution to the BSP (\ref{probeq.oper})-(\ref{probover.oper}) is constructed followed by the demonstration of the existence and uniqueness of the obtained solution. 
	\subsubsection*{Construction of the Solution:} Due to the establishment in Lemma \ref{Riesz.oper} that the collection $\{X_{0}(x), X_{1k}(x), X_{2k}(x)\}$ makes up the Riesz basis for the space $L^{2}(0,1)$, it becomes possible for us to express $u(t,x)$ and $f(x)$ in the subsequent manner
	\begin{align}
		u(t,x)&=u_{0}(t)X_{0}(x)+\sum_{k=1}^{\infty} \big(u_{1k}(t)X_{1k}(x)+u_{2k}(t)X_{2k}(x)\big),\label{u.oper}\\
		f(x)&=f_{0}X_{0}(x)+\sum_{k=1}^{\infty} \big(f_{1k}X_{1k}(x)+f_{2k}X_{2k}(x)\big).\label{f.oper}
	\end{align}
	Making use of Equations (\ref{u.oper}) and (\ref{f.oper}) in Equation (\ref{probeq.oper}) and owing to the fact that the sets $\big\{X_{0}(x), X_{1k}(x), X_{2k}(x)\big\}$ and $\big\{Y_{0}(x), Y_{1k}(x), Y_{2k}(x)\big\}$ form a bi-orthogonal system of functions for the space $L^{2}(0,1)$, the subsequent set of fractional differential equations is derived
	\begin{align}
		\mathcal{D}^{\varrho_{m}}_{t,0+}u_{0}(t)&=f_{0},\label{fdeu0.oper}\\
		\mathcal{D}^{\varrho_{m}}_{t,0+}u_{1k}(t)+\lambda_{k}u_{1k}(t)&=f_{1k},\label{fdeu1k.oper}\\
		\mathcal{D}^{\varrho_{m}}_{t,0+}u_{2k}(t)-2\sqrt{\lambda_{k}}u_{1k}(t)+\lambda_{k}u_{2k}(t)&=f_{2k}.\label{fdeu2k.oper}
	\end{align}
	Using Theorem \ref{inverse.oper} in (\ref{fdeu0.oper})-(\ref{fdeu2k.oper}), we have
	\begin{align}
		S_{\alpha_{0}+\alpha_{1}-1}\circ u_{0}(t)-S_{\alpha_{1}}\circ \phi_{0}&=f_{0},\label{invuo.oper}\\
		(S_{\alpha_{0}+\alpha_{1}-1}+\lambda_{k})\circ u_{1k}(t)-S_{\alpha_{1}}\circ \phi_{1k}&=f_{1k},\label{invu1k.oper}\\
		(S_{\alpha_{0}+\alpha_{1}-1}+\lambda_{k})\circ u_{2k}(t)-2\sqrt{\lambda_{k}}\circ u_{1k}(t)-S_{\alpha_{1}}\circ \phi_{2k}&=f_{2k}.\label{invu2k.oper}
	\end{align}
	Using Lemma \ref{Luchko-Caputo-Lemma} in (\ref{invuo.oper})-(\ref{invu2k.oper}), we can write
	\begin{align}
		u_{0}(t)=&\frac{t^{\alpha_{0}-1}}{\Gamma(\alpha_{0})}\varphi_{0}+\frac{t^{\alpha_{0}+\alpha_{1}-1}}{\Gamma(\alpha_{0}+\alpha_{1})}f_{0},\label{spaceu0.oper}\\
		u_{1k}(t)=&e_{\alpha_{0}+\alpha_{1}-1,\alpha_{0}}(t,\lambda_{k})\varphi_{1k}+e_{\alpha_{0}+\alpha_{1}-1,\alpha_{0}+\alpha_{1}}(t,\lambda_{k})f_{1k}\;,\label{spaceu1k.oper}\\
		u_{2k}(t)=&e_{\alpha_{0}+\alpha_{1}-1,\alpha_{0}}(t,\lambda_{k})\varphi_{2k}+2\sqrt{\lambda_{k}}e_{\alpha_{0}+\alpha_{1}-1,\alpha_{0}+\alpha_{1}-1}(t,\lambda_{k})*u_{1k}(t)x+e_{\alpha_{0}+\alpha_{1}-1,\alpha_{0}+\alpha_{1}}(t,\lambda_{k})f_{2k}.\label{spaceu2k.oper}
	\end{align}
	
	Taking into consideration the overdetermined condition $u(T,x)=\psi(x)$,
	\begin{align}
		f_{0}&=\frac{\Gamma(\alpha_{0}+\alpha_{1})}{T^{\alpha_{0}+\alpha_{1}-1}}\Big(\psi_{0}-\varphi_{0}\frac{T^{\alpha_{0}-1}}{\Gamma(\alpha_{0})}\Big),\label{spacef0.oper}\\
		f_{2k-1}&=	\frac{\psi_{2k-1}-\varphi_{2k-1}\;e_{\alpha_{0}+\alpha_{1}-1,\alpha_{0}}(T,\lambda_{k})}{e_{\alpha_{0}+\alpha_{1}-1,\alpha_{0}+\alpha_{1}}(T,\lambda_{k})},\label{spacef1k.oper}\\
		f_{2k}&=	\frac{\psi_{2k}-2\sqrt{\lambda_{k}}e_{\alpha_{0}+\alpha_{1}-1,\alpha_{0}+\alpha_{1}-1}(T,\lambda_{k})*u_{1k}(T)-\varphi_{2k}\;e_{\alpha_{0}+\alpha_{1}-1,\alpha_{0}}(T,\lambda_{k})}{e_{\alpha_{0}+\alpha_{1}-1,\alpha_{0}+\alpha_{1}}(T,\lambda_{k})}.\label{spacef2k.oper}
	\end{align}
	\subsubsection*{Existence of the Solution:} To establish the existence of the regular solution, we will demonstrate the continuity of the series representing $t^{\alpha_{1}}u(t,x)$, $t^{\alpha_{1}}\mathcal{D}^{\varrho_{m}}u(t,x)$, $t^{\alpha_{1}}u_{xx}(t,x)$ and $f(x)$.
	
	On using Equation (\ref{spacef0.oper}), the Cauchy-Bunyakovsky-Schwarz inequality (CBSI) and $Y_{0}(x)=1$, we have
	\begin{align}
		\big|f_{0}\big|\leq \frac{\Gamma (\alpha_{0}+\alpha_{1})}{T^{\alpha_{0}+\alpha_{1}-1}}\Big(\|\psi\|+\frac{T^{\alpha_{0}-1}}{\Gamma(\alpha_{0})}\|\phi\|\Big),\label{estf0.oper}
	\end{align}
	Due to Lemmas \ref{podlem}, \ref{AliMittagtypeLemma1}, \ref{AliMittagtypeLemma1'} and Equation (\ref{spacef1k.oper}), we obtain
	\begin{align*}
		\big|f_{1k}\big|\leq \frac{C_{1}C_{2}|\phi_{1k}|}{T^{\alpha_{1}}}+C_{2}\lambda_{k}\big|\psi_{1k}\big|.
	\end{align*}
	On using the Lemma \ref{estlemma.oper}, we get
	\begin{align}
		\big|f_{1k}\big|\leq \frac{C_{1}C_{2}\|\phi^{(2)}\|}{k^{2}T^{\alpha_{1}}}+\frac{C_{2}\|\psi^{(4)}\|}{k^{2}}.\label{estf1k.oper}
	\end{align}
	Likewise, from (\ref{spacef2k.oper}), we have
	\begin{align}
		\big|f_{2k}\big|\leq \frac{3C_{1}C_{2}\|\phi^{(2)}\|}{k^{2}T^{\alpha_{1}}}+\frac{5C_{2}\|\psi^{(4)}\|}{k^{2}}+\frac{2C_{1}\|\psi^{(3)}\|}{k^{2}}.\label{estf2k.oper}
	\end{align}
	By means of Equations (\ref{f.oper}) and (\ref{estf0.oper}) -- (\ref{estf2k.oper}), sum of the series dominating $f(x)$ has a finite upper bound. Hence, according to the WMT, $f(x)$ denotes a function that is continuous.
	
	\noindent Considering the CBSI along with (\ref{estf0.oper}) and taking into account $Y_{0}(x) = 1$, we obtain
	\begin{align}\label{spaceu0'.oper}
		t^{\alpha_{1}}\big|u_{0}(t)\big|\leq \frac{t^{\alpha_{0}+\alpha_{1}-1}\|\phi\|}{\Gamma(\alpha_{0})}+t^{\alpha_{1}}\|\psi\|+\frac{t^{\alpha_{1}}T^{\alpha_{0}-1}\|\phi\|}{\Gamma(\alpha_{0})}.
	\end{align}
	On using Lemma \ref{podlem}, Equation (\ref{spaceu1k.oper}) and (\ref{estf1k.oper}), we have
	\begin{align}
		\big|u_{1k}(t)\big|\leq 	\frac{C_{1}\big|\phi_{1k}\big|}{\lambda_{k}t^{\alpha_{1}}}+\frac{C_{1}\big|f_{1k}\big|}{\lambda_{k}}.
	\end{align}
	Using the Lemma \ref{estlemma.oper}, CBSI and the property $\big|\frac{t}{T}\big|\leq 1$, we have
	\begin{align}\label{spaceu1k'.oper}
		t^{\alpha_{1}}\big|u_{1k}(t)\big|\leq 	\frac{C_{1}\|\phi\|}{k^{2}}+\frac{C_{1}C_{2}t^{\alpha_{1}}\|\psi^{(4)}\|}{k^{4}}+\frac{C^{2}_{1}C_{2}\|\phi^{(2)}\|}{k^{4}}.
	\end{align}
	Likewise due to the Lemma \ref{podlem}, Equation (\ref{spaceu2k.oper}), (\ref{estf2k.oper}) and (\ref{spaceu1k'.oper}), we have
	\begin{align}
		t^{\alpha_{1}}\big|u_{2k}(t)\big|\leq& 	\frac{C_{1}\|\phi\|}{k^{2}}+\frac{2C^{2}_{1}C_{2}t^{\alpha_{1}}\|\psi^{(4)}\|}{k^{5}}+\frac{2C^{2}_{1}C_{2}\|\phi\|}{k^{5}}+\frac{5C_{1}C_{2}t^{\alpha_{1}}\|\psi^{(4)}\|}{k^{4}}+\frac{3C^{2}_{1}C_{2}\|\phi^{(2)}\|}{k^{4}}+\frac{2C^{2}_{1}C_{2}t^{\alpha_{1}}\|\psi^{(3)}\|}{k^{4}}.\label{spaceu2k'.oper}
	\end{align}
	By means of (\ref{spaceu0'.oper}), (\ref{spaceu1k'.oper}), and \ref{spaceu2k'.oper}) alongside the application of WMT, it is evident that $t^{\alpha_{1}}u(t,x)$ is a continuous function.
	One can observe that $D^{\alpha_{0}}_{t,0+}u(t,x)$ is bounded from above by a series that converges. Therefore, using the WMT, it can be concluded that $D^{\alpha}_{t,0+}u(t,x)$ denotes a function which is continuous. 
	
	Moreover, from Equations (\ref{fdeu0.oper})-(\ref{fdeu2k.oper}), we have the following estimates
	\begin{align}
		t^{\alpha_{1}}\big|\mathcal{D}^{\varrho_{m}}_{t,0+}u_{0}(t)\big|\leq& \frac{\Gamma(\alpha_{0}+\alpha_{1})\|\phi\|}{\Gamma(\alpha_{0})}+\frac{\Gamma(\alpha_{0}+\alpha_{1})\|\psi\|}{T^{\alpha_{0}-1}},\label{estsumu0.oper}\\
		t^{\alpha_{1}}\big|\mathcal{D}^{\varrho_{m}}_{t,0+}u_{1k}(t)\big|\leq& \sum_{k=1}^{\infty}\Big(\frac{C_{1}\|\phi^{(2)}\|}{k^{2}}+\frac{C_{1}C_{2}t^{\alpha_{1}}\|\psi^{(4)}\|}{k^{2}}+\frac{C^{2}_{1}C_{2}\|\phi^{(2)}\|}{k^{2}}+\frac{C_{2}t^{\alpha_{1}}\|\psi^{(4)}\|}{k^{2}}+\frac{C_{1}C_{2}\|\phi^{(2)}\|}{k^{2}}\Big),\label{estsumu1k.oper}\\
		t^{\alpha_{1}}\big|\mathcal{D}^{\varrho_{m}}_{t,0+}u_{2k}(t)\big|\leq& \sum_{k=1}^{\infty}\Big(\frac{3C_{1}\|\phi^{(2)}\|}{k^{2}}+\frac{4C^{2}_{1}\|\phi^{(1)}\|}{k^{2}}+\frac{2C^{2}_{1}C_{2}\|\psi^{(4)}\|}{k^{3}}+\frac{2C^{3}_{1}C_{2}\|\phi^{(2)}\|}{k^{3}}+\frac{5C_{1}C_{2}\|\psi^{(4)}\|}{k^{2}}+\frac{3C^{2}_{1}C_{2}\|\phi\|}{k^{2}}\notag\\
		&+\frac{2C^{2}_{1}C_{2}\|\psi^{(3)}\|}{k^{3}}+\frac{2C_{1}\|\phi^{(1)}\|}{k^{2}}+\frac{2C^{2}_{1}C_{2}\|\psi^{(4)}\|}{k^{2}}+\frac{2C^{2}_{1}C_{2}\|\phi^{(2)}\|}{k^{2}}+\frac{5C_{2}\|\psi^{(4)}\|}{k^{2}}+\frac{3C_{1}C_{2}\|\phi^{(2)}\|}{k^{2}}\notag\\
		&+\frac{2C_{1}C_{2}\|\psi^{(2)}\|}{k^{2}}\Big).\label{estsumu2k.oper}
	\end{align}
	By virtue of the Lemma \ref{DzhrSamkolemma} i.e. $\mathcal{D}^{\varrho_{m}}_{t,0+}\sum_{k=1}^{\infty}g_{i}(t)=\sum_{k=1}^{\infty}\mathcal{D}^{\varrho_{m}}_{t,0+}g_{i}(t)$  and (\ref{estsumu0.oper})-(\ref{estsumu2k.oper}), we notice that $t^{\alpha_{1}}\mathcal{D}^{\varrho_{m}}_{t,0+}(t,x)$ is bounded above by convergent series. Therefore, by WMT $t^{\alpha_{1}}\mathcal{D}^{\varrho_{m}}_{t,0+}(t,x)$ represents a continuous function. 
	
	\noindent Similarly for $u_{xx}(t,x)$, we have
	\begin{align*}
		\big|u_{xx}(t,x)\big|\leq& \sum_{k=1}^{\infty}\Bigg[\Bigg(\frac{C_{1}\|\phi\|}{k^{2}}+\frac{C_{1}C_{2}t^{\alpha_{1}}\|\psi^{(4)}\|}{k^{4}}+\frac{C^{2}_{1}C_{2}\|\phi^{(2)}\|}{k^{4}}\Big)+16\Big(\frac{C_{1}\|\phi\|}{k^{2}}+\frac{2C^{2}_{1}C_{2}t^{\alpha_{1}}\|\psi^{(4)}\|}{k^{5}}+\frac{2C^{2}_{1}C_{2}\|\phi\|}{k^{5}}\notag\\
		&+\frac{5C_{1}C_{2}t^{\alpha_{1}}\|\psi^{(4)}\|}{k^{4}}+\frac{3C^{2}_{1}C_{2}\|\phi^{(2)}\|}{k^{4}}+\frac{2C^{2}_{1}C_{2}t^{\alpha_{1}}\|\psi^{(3)}\|}{k^{4}}\Bigg)\Bigg],
	\end{align*}
	which on employing WMT denotes a continuous function.
	\subsubsection*{Uniqueness of the Solution:} 
	In order to demonstrate the uniqueness of $\{u(t,x), f(x)\}$, we consider two separate solutions of the system (\ref{probeq.oper})-(\ref{probover.oper}), assuming they are distinct. Following the similar approach as outlined in Theorem \ref{fp.oper}, we reach the required result.
	\end{proof}
\begin{remark}
	The recovery of the results of Kirane et al. \cite{Kirane-Salman-APC} and Furati et al. \cite{Furati-Tyiola-Kirane} can be achieved from BSP (\ref{probeq.oper})-(\ref{probover.oper}) by utilizing $\alpha_{0}=1$, $\alpha_{1}=\alpha$ and $\alpha_{0}=1-(1-\alpha)(1-\beta)$, $\alpha_{1}=1-\beta(1-\alpha)$, where $\alpha_{0},\: \alpha_{1} \in (0,1)$, respectively.
\end{remark}
\section{Concluding Remarks}\label{sec: Conclusion}
In conclusion, this article has provided a comprehensive exploration of Mikusi\'nski's type operational calculus for the Dzherbashian-Nersesian operator and its application in solving forward and backward problems. Through the utilization of this operational calculus, we derive solutions for FDEs which is expressed as infinite series of Mittag-Leffler type functions. The significance of this article lies in its potential to address complex problems in various fields, including mathematical physics and engineering, where the FDEs are encountered. Furthermore, by imposing certain conditions of the consistency and regularity of the provided data, we employed the WMT to establish the existence of a regular solution. The rigorous analysis adds a layer of confidence to the applicability of the derived solution in the practical scenarios. In summary, the operational calculus for the Dzherbashian-Nersesian operator, coupled with the use of Mittag-Leffler function and the WMT, offers a powerful toolset for solving challenging problems of fractional calculus and it's potential applications, paving the way for future advancements.

\end{document}